\documentclass[10pt,english]{amsart}
\usepackage{amsmath, amsthm, latexsym, amssymb}
\usepackage{amsfonts}
\usepackage{xypic}
\usepackage{epsfig}
\input xy
\xyoption{all}

\newcommand{\shrinkmargins}[1]{
  \addtolength{\textheight}{#1\topmargin}
  \addtolength{\textheight}{#1\topmargin}
  \addtolength{\textwidth}{#1\oddsidemargin}
  \addtolength{\textwidth}{#1\evensidemargin}
  \addtolength{\topmargin}{-#1\topmargin}
  \addtolength{\oddsidemargin}{-#1\oddsidemargin}
 \addtolength{\evensidemargin}{-#1\evensidemargin}
  }

\shrinkmargins{.7}

\theoremstyle{plain}

\newtheorem{theorem}{Theorem}[section]
\newtheorem{corollary}[theorem]{Corollary}

\newtheorem{proposition}[theorem]{Proposition}
\newtheorem{question}[theorem]{Question}

\newtheorem{definition}[theorem]{Definition}
\newtheorem{conjecture}[theorem]{Conjecture}

\theoremstyle{remark}
\newtheorem{remark}[theorem]{Remark}

\theoremstyle{definition}
\newtheorem{example}[theorem]{Example}

\def \Z { \mathbb{Z}}
\def \Q { \mathbb{Q}}

\def \R { \mathbb{R}}

\def \tr { \text{tr}}

\begin{document}
\thispagestyle{empty}
\setcounter{tocdepth}{1}

\title{On number fields with equivalent integral trace forms}
\author{Guillermo Mantilla-Soler}
\address{Department of Mathematics, University of British Columbia, 1984 Mathematics Road,\newline 
Vancouver, BC V6T 1Z2 Canada}

\email{mantilla@math.ubc.ca}

\date{}

\begin{abstract}
Let $K$ be a number field. The \textit{integral trace form} is the integral quadratic form given by $\tr_{K/\Q}(x^2)|_{O_{K}}.$ 
In this article we study the existence of non-conjugated number fields with equivalent integral trace forms. As a corollary of 
one of the main results of this paper, we show that any two non-totally real number fields with the same signature and same 
prime discriminant have equivalent integral trace forms. Additionally, based on previous results obtained by the author and the 
evidence presented here, we conjecture that any two totally real quartic fields of fundamental discriminant have equivalent trace 
zero forms if and only if they are conjugated.\\
\end{abstract}

\maketitle

\section{Introduction}

Let $K$ be a number field and let $O_{K}$ be its maximal order. 
The \textit{integral trace form} $\mathrm{tr}_{K/\Q}$, is the integral quadratic form associated to the symmetric $\Z$-bilinear form obtained by restricting the trace pairing from $K \times K$ to
$O_{K} \times O_{K}$; that is, 
\begin{displaymath}
\begin{array}{cccc}
 & O_{K} \times O_{K} &\rightarrow& \Z  \\  & (x,y) & \mapsto &
\mathrm{tr}_{K/\Q}(xy).
\end{array}
\end{displaymath}

One of the most useful invariants of a number field is its discriminant, and the significance of it can perhaps be explained in two different ways.    
On one hand, the discriminant of a number field determines which primes are ramified.  On the other hand, Hermite's Theorem shows that 
for every integer $n>0$ and for every $X>0$ 
\[ N_{n}(X):=|\{K : \mathrm{deg}(K)=n, |\mathrm{Disc(K)}| \leq X\} |< \infty.\] 
Hence we can count number fields by ordering them with respect to their discriminant,  and in fact this is how number fields are normally enumerated.  Moreover, the asymptotic behavior of $N_{n}(X)$ and related quantities have been, and still are, a very active research topic (see \cite{bha1}, \cite{bha2}, \cite{bha3}, \cite{davenport}, \cite{davenport2}, \cite{EV1}, \cite{EV2}, \cite{malle}, \cite{tu}, \cite{tho}). Since the discriminant of a number field is an invariant of the equivalence class of the integral trace form, the following question naturally arises in the study of number fields.

\begin{question}\label{Q}
Are there examples of non-conjugated number fields $K$ and $L$ such
that their corresponding integral trace forms are equivalent?
\end{question} 
  
A systematic analysis of Question \ref{Q} for cubic fields can be found in \cite{Manti}. However, at least in the generic case, not much is known about generic higher degree number fields.\footnote{In \cite{conner} the authors show that any two prime degree Galois extensions of the same degree and discriminant have equivalent trace forms.}In the next section we briefly recall certain similarities between isospectral manifolds and arithmetically equivalent number fields.\footnote{We give the corresponding definitions in the next section.} Following the above analogy between manifolds and number fields the above question could be interpreted  in Kac's language as follows: (see \S \ref{isospectral}  and \cite{Kac} for the terminology)
\[
\text{\bf Can the trace hear the shape of its field?}
\]
  
The main results of this article are the following:

\begin{theorem}\label{general}
Let $F$ and $L$ be two number fields with the same discriminant. Suppose that there is only one finite prime $p$ ramifying on $F$ and that $p$ is tamely ramified. Moreover, assume that $F$ and $L$ have the same signature, same degree and that they are non-totally real. 
Then $F$ and $L$ have equivalent integral trace forms.   
\end{theorem}

Notice that any of pair of fields $K$ and $L$ satisfying the hypotheses of Theorem \ref{general} will give a positive answer  to $\ref{Q}$. Since such examples exist (see Examples \ref{1ex}, \ref{2ex} and \ref{3ex})
we see that the trace cannot hear the shape of its field.

The \textit{trace zero module} of a number field $K$, denoted by $O_{K}^{0}$, is defined to be the set 
\[
O_{K}^{0}:=\{x \in O_K ~|~\mathrm{\tr} _{K/\Q}(x)=0\}.
\]

\begin{theorem}\label{hastadiez}
Let $n$ be a positive integer less than $11$, and let $X_n$ be the quantity described in Table $1$ below. 
Suppose that $K$ is a totally real number field of degree $n$ with fundamental
discriminant bounded by $X_n$. If $L$ is a tamely ramified number field such that there exists an
isomorphism of quadratic modules
\begin{eqnarray*}
  \langle O^{0}_{K} , \mathrm{tr}_{K/ \mathbb{Q}}(x^2)|_{O^{0}_{K}} \rangle & \cong & \langle
O^{0}_{L},\mathrm{tr}_{L/ \mathbb{Q}}(x^2)|_{O^{0}_{L}} \rangle,
\end{eqnarray*}
then $K \cong L$.

\begin{table}[h!]
\begin{center}
{
\begin{tabular}[t]{|c|c|}
	\hline
  $ X_{n}$ & $ n $  \\
	\hline \hline
$\infty$ & $ 1,2,3$   \\
	\hline	
$ 1 \times 10^9$ &  $ 4,5,6$   \\
	\hline
$8.9\times 10^{10}$ & $ 7$   \\ 
	\hline
$2.5\times 10^9$ & $ 8$     \\
        \hline
$2.8\times 10^{10}$  & $ 9$    \\
        \hline    
$2.8\times 10^{11}$ & $10$       \\
        \hline              
\end{tabular}
}
\end{center}
\caption{\quad}
\end{table}

\end{theorem}

\section{Overview}
In this section we explain how Question \ref{Q} leads naturally to Theorems \ref{general} and \ref{hastadiez} and we state some of their consequences.  In addition, as a result of our study of Question \ref{Q}, we observe some interesting analogies between number fields and certain types of manifolds.\\ 

One difficulty that arises in looking for answers to Question \ref{Q} is that in principle it is hard to find non-conjugated number fields sharing the same discriminant. Via class field 
theory  one can determine when a fundamental discriminant $d$ is realized as the discriminant of more than one cubic field. In fact, one can 
produce similar arguments to decide the same about quartic number fields. However, since $S_{n}$ is not solvable for $n>4$, this method is not useful anymore. Nevertheless, there is another way to find number fields  sharing discriminants which a priori works in any degree. 

\begin{definition}
Two number fields $K$ and $L$ are called \textit{arithmetically
equivalent}\footnote{The definition given here, although equivalent, is not what was initially called arithmetical equivalence. See \cite{Perlis1} for the orignial definition.} if their Dedekind zeta functions $\zeta_{K}$ and $\zeta_{L}$ 
coincide.
\end{definition}
By a well known result of Perlis arithmetically equivalent number fields have the same degree and discriminant. 
Therefore, the following theorem of Perlis (\cite{Perlis1}) gives a simple
tool to find candidates to test Question \ref{Q}.

\begin{theorem}[Perlis]\label{quasiconjugated}
Let $K$ and $L$ be number fields and let $N$ to be the Galois closure of
$KL$ over $\Q$. Let $G =\mathrm{Gal}(N/\Q)$, $H=\mathrm{Gal}(N/K)$ and
$I=\mathrm{Gal}(N/L)$. Then $K$ and $L$ are arithmetically equivalent if
and only if $H$ and $I$ are almost conjugated\footnote{ Recall that a pair of subgroups 
$H$ and $I$ of a finite group $G$ are \emph{almost conjugate} if and only 
if $\mathrm{Ind}_{H}^{G}1_{H}=\mathrm{Ind}_{I}^{G}1_{I}.$} subgroups of $G$.
\end{theorem}

It follows that we can find non-conjugated number fields of the same degree and discriminant by finding groups $G$ that are verifiable Galois over $\Q$
containing a pair of non-conjugated, almost conjugated subgroups. For number fields of small degree, namely less than $11$, the group $G$ of
smallest order with this property is $G=\Z/8\Z \rtimes (\Z/8\Z)^{*}$ (see \cite[Chapter III, Theorem 1.14]{Klingen}).  
In Example \ref{trazaocho} we exhibit, based on a technique developed by Gerst and Schinzel \cite{Gertz}, two non-conjugated, almost 
conjugated octic number fields of discriminant $2^{10}\cdot3^7\cdot5^7,$ such that the Galois closure of their compositum has Galois group 
isomorphic to $G$. These number fields have been chosen so that the trace form of one reduces to the zero form module $2$ while the other does not. In particular it is trivial to see that their integral trace forms are not equivalent.

\begin{example}\label{trazaocho}
Let $F$ and $L$ be the number fields defined by the polynomials
$p_F=x^8+15$ and $p_L:=x^8+240$, respectively. Their common Galois closure over $\Q$ is the field $N=F(\zeta_{8})$ which has Galois group the semidirect
product $G=\Z/8\Z \rtimes (\Z/8\Z)^{*}$. Since $G$ can be described as the group of affine transformations of $\Z/8\Z$, it can be shown, in the exact same way as in \cite[Chapter III, \S1 Example 1.9 ]{Klingen}, that $F$ and $L$ are almost conjugated, and thus arithmetically equivalent.  
Let $y$ be a root of $p_F$ and $z$ a root $p_L$. The following are integral
bases for $F$ and $L$ respectively:

\begin{eqnarray*}
 B_F &=&\left \{ 1,
    y,
    y^2,
    y^3,
    \frac{1}{2}(y^4 + 1),
    \frac{1}{2}(y^5 + y),
    \frac{1}{4}(y^6 + y^4 + y^2 + 1), \right. \\ & & 
   \left. \frac{1}{8}(y^7 + y^6 + y^5 + y^4 + y^3 + y^2 + y + 1) \right \}\\ 
B_L & = & \left \{1,
    z,
    \frac{1}{2}z^2,
    \frac{1}{4}(z^3 + 2z),
    \frac{1}{8}(z^4 + 4),
    \frac{1}{16}(z^5 + 4z^2 + 12z + 8),
    \frac{1}{32}(z^6 + 2z^4 + 4z^2 + 8), \right. \\ & &
   \left. \frac{1}{64}(z^7 + 2z^5 + 4z^4 + 12z^3 + 16z^2 + 24z + 16) \right \}.
\end{eqnarray*} The Gramm matrices of the trace form in these bases are the following: 
\begin{footnotesize}
\[M_F =
\begin{pmatrix}
  8 & 0 & 0 & 0 & 4 & 0 & 2 & 1 \\
  0 & 0 & 0 & 0 & 0 & 0 & 0 & -15 \\
  0 & 0 & 0 & 0 & 0 & 0 &-30& -15 \\
  0 & 0 & 0 & 0 & 0 &-60& 0 & -15 \\
  4 & 0 & 0 & 0 &-28& 0 &-14& -7 \\
  0 & 0 & 0 &-60& 0 & 0 & 0 & -15 \\
  2 & 0 &-30& 0 &-14& 0 &-22& -11 \\
  1 &-15&-15&-15& -7&-15&-11& -13 
\end{pmatrix} M_L  = 
\begin{pmatrix}
  8 & 0 & 0 & 0 & 4 & 4 & 2 & 2 \\
  0 & 0 & 0 & 0 & 0 & 0 & 0 & -30 \\
  0 & 0 & 0 & 0 & 0 & 0 &-30& 0  \\
  0 & 0 & 0 & 0 & 0 &-30& 0 & -30 \\
  4 & 0 & 0 & 0 &-28& 2 &-14& -14 \\
  4 & 0 & 0 &-30& 2 & 2 &-14& -44 \\
  2 & 0 &-30& 0 &-14&-14&-22& -22 \\
  2 &-30& 0 &-30&-14&-44&-22& -52 
\end{pmatrix}\]
\end{footnotesize}
\noindent Since $M_L$  clearly is zero as an element of 
M$_{2}(\Z/2\Z)$ 
and $M_{F}$ is not,  the integral trace forms of $F$ and $L$ are not equivalent. \end{example}

\subsection{Isospectral manifolds}\label{isospectral} In this subsection we explore the connections between number fields with the same discriminant and 
isospectral manifolds.\\

Let $M$ a closed Riemannian manifold of dimension $n$ and associated Laplacian $\Delta$. Note that the eigenvalue spectrum of 
$\Delta$ acting on $L_{2}(M)$ is discrete with each eigenvalue occurring with finite multiplicity.   

\begin{definition}
Suppose that $M_{1}$ and $M_{2}$ are two closed Riemannian $n$-manifolds. Then $M_1$ and $M_2$ are said to be \textit{isospectral} if they 
have equal eigenvalue spectra for their respective Laplacians.
\end{definition}

Although isospectral manifolds and number fields of the same discriminant are, as far as we know, unrelated, it turns out that finding number fields with the same discriminant by the method described above is analogous to Sunada's method of finding isospectral Riemannian manifolds.

The first example of two non-isometric isospectral manifolds was given by Milnor \cite{Milnor} in dimension $16$. Later Sunada \cite{Su} produced a general method for constructing isospectral manifolds. In fact, Milnor's construction 
can be seen as a particular case of Sunada's method. Nevertheless, the problem of finding two non-isometric, isospectral Riemannian $2$-manifolds  remained open.

Roughly speaking, isospectral $2$-manifolds are drums that sound the same (see \cite{Kac}). Finally the matter was settled by Gordon, Webb and Wolpert \cite{GWW} when they answered Kac's famous question ``Can one hear the shape of a drum?"\footnote{Kac asked in \cite{Kac} if it is possible to have two differently shaped drums which have the same eigenspectrum.} by constructing, based on Sunada's method, a pair of regions in the plane that have different shapes but identical eigenvalues.

\begin{theorem}[Sunada's method]
Let $M$ be a closed Riemannian manifold
such that there exists a surjective homomorphism \[\phi : \pi_1(M) \rightarrow G.\]
Suppose that $H$ and $I$ are almost conjugate subgroups of $G$.  If $M_H$ and $M_I$ are the 
finite covers of $M$ associated to the finite index subgroups $\Phi^{-1}(H)$ and $\Phi^{-1}(I)$, 
then $M_H$ and $M_I$ are isospectral.
\end{theorem}

It follows that Sunada's method, together with Theorem \ref{quasiconjugated}, yields to an analogy between isospectral manifolds and arithmetically equivalent number fields. Hence, the following question arises naturally. 
\begin{question}
Do two non-conjugated, arithmetically 
equivalent number fields $F$ and $L$ exist such that their corresponding trace forms are equivalent? 
\end{question}
Notice that a positive answer to this question will imply a positive answer to Question \ref{Q}. 

The first case in which the above question is interesting is for septic number fields since any two arithmetically equivalent number fields of degree less than $7$ are conjugated (see for example \cite[Chapter III, Theorem 1.16]{Klingen}).  
 
\begin{proposition}\label{siete}
There exist two arithmetically equivalent, non-conjugated, septic number fields such that their integral trace forms are equivalent.
\end{proposition}

In other words, even under arithmetically equivalence the trace cannot hear the shape of its field.

\subsection{The right question}
As seen in \cite{Manti} when dealing with Question \ref{Q} it is convenient to separate the problem into two different cases, namely, the positive and negative discriminant cases. If instead of the sign of the discriminant we consider the definiteness of the trace, then we get a more suitable division of the problem, which
agrees with the sign of the discriminant in the cubic case.  The following proposition shows that splitting the problem between fields with positive definite trace and non-positive definite trace is equivalent to dividing the problem between fields that are  totally real and those that are not. \

\begin{proposition}[Taussky]\label{signature}
For every number field $F$ the signature of $\mathrm{tr}_{F/\Q}$ is equal to the number of real places $r_{F}$ of $F$. 
\end{proposition}

Let us first consider the case of number fields in which the trace form is non-definite, or equivalently by Proposition \ref{signature}, number fields that are non-totally real. In this case we ask: 

Are there examples of  non-conjugated,  non-totally real number fields $F$ and $L$  such
that their corresponding trace forms are equivalent? 
Moreover, can we get them to have fundamental discriminant?  
In \cite{Manti} this question 
is answered affirmatively for
cubic fields using the trace zero form. Theorem \ref{general} provides a more general answer  to the above questions, and in fact  some of the examples given in 
\cite{Manti} can also be deduced from Theorem \ref{general}.\\

Since examples of pairs of fields $F$ and $L$ satisfying the hypotheses of Theorem \ref{general} exist (see Examples \ref{1ex}, \ref{2ex} and \ref{3ex}), we conclude that for non-totally real number fields the trace form does not characterize its number field. In view of this result, and the results proven in \cite{Manti}, the following question is relevant:

\begin{question}\label{pregunta}
Are the examples of non-conjugated totally real number fields $F$ and $L$ of fundamental discriminant such
that their corresponding integral trace zero forms are equivalent?
\end{question}

Theorem \ref{hastadiez} gives a partial answer to Question \ref{pregunta} in the sense that for $n <11$, and for the restricted discriminants given in the hypothesis, the answer to Question \ref{pregunta} is no.
We remark that the implicit result that Theorem \ref{hastadiez} seems to imply, which is obviously true for $n=1,2$  and already known for $n=3$, should be considered mainly for the case $n=4.$  To understand this point it is important to notice that the data for $n \ge 5$ is very scarce. For $5 \leq n \leq 10$, and the values of $X_n$ stated in Theorem \ref{hastadiez}, there are $5,122$ number fields, all of them for $n \in \{5,6,7\}$, with non-trivial discriminant multiplicity.\footnote{A number field $F$ is said to have \emph{trivial discriminant multiplicity} if for all number field $K$ of the same signature and Galois group as $F$ one has that $\mathrm{Disc}(F)=\mathrm{Disc}(K)$ implies that $F$ is conjugated to $K$.} In contrast, among all of the totally real quartic number fields of discriminant up to $10^9$, there are  $1,301,472$ number fields with non-trivial discriminant multiplicity. Therefore, Theorem \ref{hastadiez} can be thought as numerical evidence of the following:\\

\begin{conjecture} 
Let $K$ be a totally real quartic number field with fundamental
discriminant. If $L$ is a tamely ramified number field such that an
isomorphism of quadratic modules
\begin{eqnarray*}
  \langle O^{0}_{K} , \tr_{K/ \mathbb{Q}}(x^2)|_{O^{0}_{K}} \rangle & \cong & \langle
O^{0}_{L},\tr_{L/ \mathbb{Q}}(x^2)|_{O^{0}_{L}} \rangle,
\end{eqnarray*}
exists, then $K \cong L$.
\end{conjecture}

\section{Proofs and examples}
We now proceed to prove all the claims made in the previous sections.
\subsubsection{Proof of Proposition \ref{siete} }
\begin{proof}
Let $F$ and $L$ be the number fields defined by the polynomials
\begin{align*}
p_F&=x^7-3x^6+4x^4+x^3-4x^2-x+1,\\
p_L&=x^7-3x^6+2x^5+4x^4-3x^3-2x^2-x-1,
\end{align*}
respectively.  Since the factorization type of $2741$ on $O_F$ is $(1,1,1,2)$ 
while its factorization type on $O_L$ is $(2,2,1,1)$,  we see that $F$ and $L$ are not conjugated.  A calculation shows that $[FL:\Q] \leq 28$, thus 
$F$ and $L$ are not linearly disjoint. Since $F$ and $L$ have prime degree over $\Q$ it follows from \cite{Perlis2} that $F$ and $L$ are arithmetically 
equivalent.  Finally, since $F$ and $L$ have both signature $(3,2)$ and discriminant $(2741)^2$, we can apply Theorem \ref{general} to conclude that $F$ 
and $L$ have equivalent integral trace forms. 
\end{proof}

\subsubsection{Proof of Theorem \ref{general}.} 
\begin{proof}  We may assume that 
$[F:\Q]>2$. Let tr$_{F/\Q}$ and tr$_{L/\Q}$ be the rational quadratic forms given by the trace forms over the fields $F$ and $L$.  First 
we show that these two rational quadratic forms are rationally  equivalent. It follows from the signature hypothesis that tr$_{F/\Q}$ and tr$_{L/\Q}$  have the same dimension. Hence  tr$_{F/\Q}$ and tr$_{L/\Q}$ are rationally equivalent if and only if their Hasse invariant $\mathrm{h}_{\ell}$ is the same for every prime $\ell$, finite or infinite. We now show that $\mathrm{h}_{\ell}(\mathrm{tr}_{F/\Q})=\mathrm{h}_{\ell}(\mathrm{tr}_{L/\Q})$ for all possible values of $\ell$.
\begin{itemize}
\item $\ell= \infty.$  If we denote the signature of $F$ by $(r_{F}, s_{F})$ then, by Proposition \ref{signature}, we see 
that $\mathrm{tr}_{F/\Q}$ is equivalent over $\R$ to a form of type $(r_{F} + s_{F}, s_{F})$. In particular, \[\mathrm{h}_{\ell}(\mathrm{tr}_{F/\Q})= (-1)^{\frac{s_{F}(s_{F}-1)}{2}}=(-1)^{\frac{s_{L}(s_{L}-1)}{2}}=\mathrm{h}_{\ell}(\mathrm{tr}_{L/\Q}).\]
\item $\ell= 2.$  For tamely ramified extensions, the determinant and the discriminant determine the trace form as a quadratic form over $\Z_{2}$ (see  \cite[remark above Lemme 1.11]{emp}) In particular, due to the ramification hypothesis,  we see that $\mathrm{h}_{2}(\mathrm{tr}_{F/\Q})=\mathrm{h}_{2}(\mathrm{tr}_{L/\Q})$. 
\item $\ell \neq p,2.$   Let $d$ be the common discriminant. Since  $\ell$ does not divide $2d$ we see that
$\mathrm{h}_{\ell}(\mathrm{tr}_{F/\Q})=1$ and $\mathrm{h}_{\ell}(\mathrm{tr}_{L/\Q})=1$. 
\item $\ell= p.$
Since $\mathrm{h}_{\ell}(\mathrm{tr}_{F/\Q})=\mathrm{h}_{\ell}(\mathrm{tr}_{L/\Q})$ for all places $\ell \neq p$, it follows from the product formula that $\mathrm{h}_{p}(\mathrm{tr}_{F/\Q})=\mathrm{h}_{p}(\mathrm{tr}_{L/\Q})$. 
\end{itemize} 
Recall that for $F/\Q$ and $L/\Q$ tamely ramified extensions, the forms  tr$_{F/\Q}$ and tr$_{L/\Q}$  are in the same 
genus if and only if the fields have equal discriminant and the forms are rationally equivalent (see \cite[Theorem 1.10]{emp}). Since we just showed that tr$_{F/\Q}$ and tr$_{L/\Q}$ are rationally equivalent we conclude that they are in the same genus. Our next step is to show that the forms are in fact in the same spinor genus.  Let $d=p^m$ be the common discriminant of $F$ and $L$. Since $F/\Q$ is tamely ramified, the valuation of $d$ at $p$ is given by $n-t$, where $n=[F:\Q]$ and $t$ is the sum of residue degrees at $p$ (see \cite[Chapter III]{Serre2}). Therefore, $m<n$ and in 
particular we see that $d$ is not divisible by the $\frac{n(n-1)}{2}$-th  power for any integer greater than $2$. Hence, it follows from 
\cite[Chapter 7, \S3, Corollary 1]{Watson} that the genus and the spinor genus of tr$_{F/\Q}$ coincide. Since the same is true for tr$_{L/\Q}$, and  since tr$_{F/\Q}$ and tr$_{L/\Q}$ are in the same genus, we see that they are in the same spinor genus. Since $3\leq n$ and tr$_{F/\Q}$, tr$_{L/\Q}$ 
are indefinite integral forms, they are equivalent.
\end{proof}

\begin{corollary}\label{exampleswithprimes}
Any two cubic fields with the same negative prime discriminant have equivalent integral trace forms. More generally, any two non-totally real number 
fields with the same signature and same prime discriminant have equivalent integral trace forms.  
\end{corollary}

\begin{example}\label{1ex}
The following polynomials describe four non-conjugated cubic fields of prime discriminant $-3299$. By Corollary \ref{exampleswithprimes} we know that they have
equivalent integral trace forms.
\begin{itemize}
\item $x^3 + 2x+11$
\item $x^3 - 16x + 27$
\item $x^3 - x^2 + 9x - 8$
\item $x^3 - x^2 + 3x + 10$
\end{itemize}
\end{example}

\begin{example}\label{2ex}
The following polynomials describe two non-conjugated quartic fields of prime discriminant $7537$. Since they are non-totally real,  by Corollary \ref{exampleswithprimes},  they have
equivalent integral trace forms.
\begin{itemize}
\item $x^4 - x^3 + 5x^2 - 4x + 5$
\item $x^4 + 4x^2 - 5x + 2$
\end{itemize}
\end{example}

\begin{example}\label{3ex}
The following polynomials describe two non-conjugated quintic fields of prime discriminant $34129$. Since they are non-totally real, by Corollary \ref{exampleswithprimes}, that they have
equivalent integral trace forms.
\begin{itemize}
\item $x^5 + 2x^3 - x^2 - 2x - 1$
\item $x^5 - 2x^4 + x^2 - 2x + 4$
\end{itemize}
\end{example}

\begin{remark}
Notice that if we remove the non-totally real condition in Theorem \ref{general}, we no longer  can conclude that the trace forms are equivalent. 
However, from the proof of the theorem we conclude that tr$_{F/\Q}$ and tr$_{L/\Q}$ are in the same spinor genus. 
\end{remark}

\begin{example}
The following polynomials define non-conjugated quartic fields such that their trace forms are pairwise non-equivalent quadratic forms in the same spinor genus. 
\begin{itemize}
\item $19x^4 - x^3 - 10x^2 + 1$
\item $16x^4 -23x^2 - 18x + 1$
\item $17x^4 +8 x^3 -10x^2 - x + 1$
\end{itemize}
\end{example}

\subsubsection{Proof of Theorem \ref{hastadiez}.}

\begin{proof}
First we show that in Theorem \ref{hastadiez} we can replace ``tamely ramified number field" with ``number field with fundamental discriminant". Let $F$ be a number field and let us denote by  $\mathrm{tr}^{0}_{F/\Q}$ its trace zero form. If $F$ is tamely ramified we see, from \cite[Lemma 1]{Mau}, that the trace map $\mathrm{tr}_{F/\Q}: O_{F} \rightarrow \Z$ is surjective. In particular, the $\Z$-submodule of $O_{F}$ generated by $\Z$ and $O^{0}_{F}$ is a submodule of index $[F:\Q]$. It follows that  $|\mathrm{disc}( \mathrm{tr}^{0}_{F/\Q})|=[F:\Q]|\mathrm{disc}( \mathrm{tr}_{F/\Q})|$ for tamely ramified $F$. Hence if $K$ is a totally real number field with fundamental discriminant and $L$ is a tamely ramified number field, an isomorphism between their respective trace zero forms not only implies that $K$ and $L$ have the same degree but also the same discriminant. In particular, $L$ must have fundamental discriminant. We can now proceed assuming that ``number field with fundamental discriminant" replaces ``tamely ramified number field". We notice that the cases $n=1,2$ are trivial and the case $n=3$ is  \cite[Theorem 6.5]{Manti}.  In the remaining cases we have written code in \cite{magma} to verify the theorem (see next section).  Gunter Malle has provided us, via personal communication (see also his online database \cite{KMa}), with tables of totally real quartic, quintic, sextic and septic number fields with non-trivial discriminant multiplicity and discriminant up to $X_n$.  We have also used data from \cite{Diaz}, \cite{Jones} and \cite{Voight}, which together with the tables from Kl\"{u}eners-Malle give a complete list of all totally real number fields of degree $n$ with discriminant up to $X_{n}$.
\end{proof}

\subsection{{\tt Magma} code }

The following is the code that we have written to verify whether or not two totally real number fields have conjugated trace zero forms.  It is important to observe that 
it is assumed that the input fields are totally real.
\begin{itemize}
\item Given a number field $F$ defined by a polynomial $f$ the function \verb"TraceZeroSubspace"  calculates a $\Z$-basis for $O_{F}^{0}.$  The input of the function is $f$. 
\begin{verbatim} 
TraceZeroSubspace := function(f);
F:=NumberField(f);
O_F:=MaximalOrder(F);
Z_Fbas :=Basis(O_F);
M := Matrix(IntegerRing(),[[Trace(x) : x in Z_Fbas]]);
kerM := Kernel(Transpose(M));
return[&+[v[i]*Z_Fbas[i] : i in[1..Degree(Z_F)]] : v in Basis(kerM)];
end function;
\end{verbatim}
\item Given a number field $F$ defined by a polynomial $f$, the function \verb"TraceMatrixZero" calculates the Gramm matrix for the trace zero form in a basis given by the function 
\verb"TraceZeroSubspace".
The input of the function is $f$.
\begin{verbatim}
TraceMatrixZero :=function(f);
Z_F:=TraceZeroSubspace(f);
a:=#Z_F;
MF:=Matrix(a, [Trace(Z_F[i]*Z_F[j]) : i,j in [1..a]]);
return MF;
end function;
\end{verbatim}
\item  Given two number fields $F,K$ defined by polynomials $f$ and $k$, the following function returns \verb"True"  if and only if the integral trace zero forms of $F$ and $K$ are equivalent, otherwise it returns \verb"False". The input of the function is the pair $f,k$.
\begin{verbatim}
IsIntegralTraceZeroIsometric :=function(f,k);
MF:=TraceMatrixZero(f);
MK:=TraceMatrixZero(k);
LF:=LatticeWithGram(MF);
LK:=LatticeWithGram(MK);
return IsIsometric(LF,LK);
end function;
\end{verbatim}
\end{itemize}

\section*{Acknowledgements}

I thank the Mathematics departments at UW-Madison and UBC-Vancouver for granting me access to their computer servers to carry out some of the computations needed in the paper.  I'm  specially thankful to Nigel Boston for allowing me to use his personal servers where the longer calculations were made. I also like to express my gratitude to Gunter Malle for kindly answering all my many requests for tables of number fields with specific conditions, and to John Voight  for his help with {\tt Magma}. Finally I'd like to thank to Jordan Ellenberg, Andrea, Jose and Samantha for their valuable comments on a draft of this paper.

\end{document}